\documentclass[12pt]{article}
\usepackage{latexsym, amssymb}

\usepackage{mathbbol}

\DeclareMathAlphabet\gothic{U}{euf}{m}{n}

\makeatletter
\def\eqnarray{\stepcounter{equation}\let\@currentlabel=\theequation
\global\@eqnswtrue
\tabskip\@centering\let\\=\@eqncr
$$\halign to \displaywidth\bgroup\hfil\global\@eqcnt\z@
  $\displaystyle\tabskip\z@{##}$&\global\@eqcnt\@ne
  \hfil$\displaystyle{{}##{}}$\hfil
  &\global\@eqcnt\tw@ $\displaystyle{##}$\hfil
  \tabskip\@centering&\llap{##}\tabskip\z@\cr}

\def\endeqnarray{\@@eqncr\egroup
      \global\advance\c@equation\m@ne$$\global\@ignoretrue}

\def\@yeqncr{\@ifnextchar [{\@xeqncr}{\@xeqncr[5pt]}}
\makeatother

\begin{document}
\bibliographystyle{tom}

\newtheorem{lemma}{Lemma}[section]
\newtheorem{thm}[lemma]{Theorem}
\newtheorem{cor}[lemma]{Corollary}
\newtheorem{voorb}[lemma]{Example}
\newtheorem{rem}[lemma]{Remark}
\newtheorem{prop}[lemma]{Proposition}
\newtheorem{stat}[lemma]{{\hspace{-5pt}}}
\newtheorem{obs}[lemma]{Observation}
\newtheorem{defin}[lemma]{Definition}

\newenvironment{remarkn}{\begin{rem} \rm}{\end{rem}}
\newenvironment{exam}{\begin{voorb} \rm}{\end{voorb}}
\newenvironment{defn}{\begin{defin} \rm}{\end{defin}}
\newenvironment{obsn}{\begin{obs} \rm}{\end{obs}}

\newcommand{\gota}{\gothic{a}}
\newcommand{\gotb}{\gothic{b}}
\newcommand{\gotc}{\gothic{c}}
\newcommand{\gote}{\gothic{e}}
\newcommand{\gotf}{\gothic{f}}
\newcommand{\gotg}{\gothic{g}}
\newcommand{\gothh}{\gothic{h}}
\newcommand{\gotk}{\gothic{k}}
\newcommand{\gotm}{\gothic{m}}
\newcommand{\gotn}{\gothic{n}}
\newcommand{\gotp}{\gothic{p}}
\newcommand{\gotq}{\gothic{q}}
\newcommand{\gotr}{\gothic{r}}
\newcommand{\gots}{\gothic{s}}
\newcommand{\gotu}{\gothic{u}}
\newcommand{\gotv}{\gothic{v}}
\newcommand{\gotw}{\gothic{w}}
\newcommand{\gotz}{\gothic{z}}
\newcommand{\gotA}{\gothic{A}}
\newcommand{\gotB}{\gothic{B}}
\newcommand{\gotG}{\gothic{G}}
\newcommand{\gotL}{\gothic{L}}
\newcommand{\gotS}{\gothic{S}}
\newcommand{\gotT}{\gothic{T}}

\newcounter{teller}
\renewcommand{\theteller}{\Roman{teller}}
\newenvironment{tabel}{\begin{list}%
{\rm \bf \Roman{teller}.\hfill}{\usecounter{teller} \leftmargin=1.1cm
\labelwidth=1.1cm \labelsep=0cm \parsep=0cm}
                      }{\end{list}}

\newcounter{tellerr}
\renewcommand{\thetellerr}{(\roman{tellerr})}
\newenvironment{subtabel}{\begin{list}%
{\rm  (\roman{tellerr})\hfill}{\usecounter{tellerr} \leftmargin=1.1cm
\labelwidth=1.1cm \labelsep=0cm \parsep=0cm}
                         }{\end{list}}
\newenvironment{ssubtabel}{\begin{list}%
{\rm  (\roman{tellerr})\hfill}{\usecounter{tellerr} \leftmargin=1.1cm
\labelwidth=1.1cm \labelsep=0cm \parsep=0cm \topsep=1.5mm}
                         }{\end{list}}

\newcommand{\Ni}{{\bf N}}
\newcommand{\Ri}{{\bf R}}
\newcommand{\Ci}{{\bf C}}
\newcommand{\Ti}{{\bf T}}
\newcommand{\Zi}{{\bf Z}}
\newcommand{\Fi}{{\bf F}}

\newcommand{\proof}{\mbox{\bf Proof} \hspace{5pt}} 
\newcommand{\remark}{\mbox{\bf Remark} \hspace{5pt}}
\newcommand{\ruimte}{\vskip10.0pt plus 4.0pt minus 6.0pt}

\newcommand{\simh}{{\stackrel{{\rm cap}}{\sim}}}
\newcommand{\ad}{{\mathop{\rm ad}}}
\newcommand{\Ad}{{\mathop{\rm Ad}}}
\newcommand{\Aut}{\mathop{\rm Aut}}
\newcommand{\arccot}{\mathop{\rm arccot}}
\newcommand{\capp}{{\mathop{\rm cap}}}
\newcommand{\rcapp}{{\mathop{\rm rcap}}}
\newcommand{\diam}{\mathop{\rm diam}}
\newcommand{\divv}{\mathop{\rm div}}
\newcommand{\codim}{\mathop{\rm codim}}
\newcommand{\RRe}{\mathop{\rm Re}}
\newcommand{\IIm}{\mathop{\rm Im}}
\newcommand{\Tr}{{\mathop{\rm Tr}}}
\newcommand{\Vol}{{\mathop{\rm Vol}}}
\newcommand{\card}{{\mathop{\rm card}}}
\newcommand{\supp}{\mathop{\rm supp}}
\newcommand{\sgn}{\mathop{\rm sgn}}
\newcommand{\essinf}{\mathop{\rm ess\,inf}}
\newcommand{\esssup}{\mathop{\rm ess\,sup}}
\newcommand{\Int}{\mathop{\rm Int}}
\newcommand{\Leibniz}{\mathop{\rm Leibniz}}
\newcommand{\lcm}{\mathop{\rm lcm}}
\newcommand{\loc}{{\rm loc}}

\newcommand{\mod}{\mathop{\rm mod}}
\newcommand{\spann}{\mathop{\rm span}}
\newcommand{\one}{\mathbb{1}}

\hyphenation{groups}
\hyphenation{unitary}

\newcommand{\tfrac}[2]{{\textstyle \frac{#1}{#2}}}

\newcommand{\cb}{{\cal B}}
\newcommand{\cc}{{\cal C}}
\newcommand{\cd}{{\cal D}}
\newcommand{\ce}{{\cal E}}
\newcommand{\cf}{{\cal F}}
\newcommand{\ch}{{\cal H}}
\newcommand{\ci}{{\cal I}}
\newcommand{\ck}{{\cal K}}
\newcommand{\cl}{{\cal L}}
\newcommand{\cm}{{\cal M}}
\newcommand{\co}{{\cal O}}
\newcommand{\cs}{{\cal S}}
\newcommand{\ct}{{\cal T}}
\newcommand{\cx}{{\cal X}}
\newcommand{\cy}{{\cal Y}}
\newcommand{\cz}{{\cal Z}}

\newcommand{\wtozp}{W^{1,2}\raisebox{10pt}[0pt][0pt]{\makebox[0pt]{\hspace{-34pt}$\scriptstyle\circ$}}}
\newlength{\hightcharacter}
\newlength{\widthcharacter}
\newcommand{\covsup}[1]{\settowidth{\widthcharacter}{$#1$}\addtolength{\widthcharacter}{-0.15em}\settoheight{\hightcharacter}{$#1$}\addtolength{\hightcharacter}{0.1ex}#1\raisebox{\hightcharacter}[0pt][0pt]{\makebox[0pt]{\hspace{-\widthcharacter}$\scriptstyle\circ$}}}
\newcommand{\cov}[1]{\settowidth{\widthcharacter}{$#1$}\addtolength{\widthcharacter}{-0.15em}\settoheight{\hightcharacter}{$#1$}\addtolength{\hightcharacter}{0.1ex}#1\raisebox{\hightcharacter}{\makebox[0pt]{\hspace{-\widthcharacter}$\scriptstyle\circ$}}}
\newcommand{\scov}[1]{\settowidth{\widthcharacter}{$#1$}\addtolength{\widthcharacter}{-0.15em}\settoheight{\hightcharacter}{$#1$}\addtolength{\hightcharacter}{0.1ex}#1\raisebox{0.7\hightcharacter}{\makebox[0pt]{\hspace{-\widthcharacter}$\scriptstyle\circ$}}}

\thispagestyle{empty}

\begin{center}
{\Large\bf Conservation and invariance properties \\[3mm]
of submarkovian semigroups} \\[5mm]
\large A.F.M. ter Elst$^1$ and  Derek W. Robinson$^2$

\end{center}

\vspace{5mm}

\begin{center}
{\bf Abstract}
\end{center}

\begin{list}{}{\leftmargin=1.8cm \rightmargin=1.8cm \listparindent=10mm 
   \parsep=0pt}
\item
Let $\ce$ be a Dirichlet form on $L_2(X)$ and $\Omega$ an open subset of $X$.
Then one can define Dirichlet forms $\ce_D$,  or $\ce_N$, corresponding to 
$\ce$ but with Dirichlet, or Neumann, boundary conditions imposed on the boundary
$\partial\Omega$ of $\Omega$.
If $S$, $S^D$ and $S^N$ are the associated submarkovian semigroups
we prove, under general assumptions of regularity and locality, that 
$S_t\varphi =S^D_t\varphi$ for all $\varphi\in L_2(\Omega)$ and 
$t>0$ if and only if the
capacity $\capp_\Omega(\partial\Omega)$ of $\partial \Omega$ relative
to $\Omega$ is zero.
Moreover, if $S$ is conservative, i.e.\ stochastically complete,
then  $\capp_\Omega(\partial\Omega)=0$ if and only if $S^D$ is 
conservative on $L_2(\Omega)$.
Under slightly more stringent assumptions we also prove that the 
vanishing of the relative capacity is equivalent to $S^D_t \varphi = S^N_t \varphi$
 for all $\varphi\in L_2(\Omega)$ and  $t>0$.

\end{list}

\vspace{7cm}
\noindent
February 2009

\vspace{5mm}

\vspace{5mm}

\noindent
{\bf Home institutions:}    \\[3mm]
\begin{tabular}{@{}cl@{\hspace{10mm}}cl}
1. & Department of Mathematics  &
  2. & Centre for Mathematics and its Applications  \\
& University of Auckland   & 
& Mathematical Sciences Institute  \\
& Private bag 92019  &
  & Australian National University  \\
  & Auckland & &
   Canberra, ACT 0200  \\
& New Zealand  & &
  Australia  \\
  \end{tabular}

\newpage
\setcounter{page}{1}

\section{Introduction}\label{Sord1}

In two earlier papers \cite{RSi} \cite{ER29} the relationships 
between the invariance of a set $\Omega$
under the action of a submarkovian semigroup $S$ and capacity conditions 
on the boundary
$\partial\Omega$ of the set were explored.
In the current paper we demonstrate  that these features are 
connected to the conservative 
property for the semigroup $S^D$ obtained by imposing 
Dirichlet boundary conditions on 
$\partial\Omega$.
Under quite general conditions the latter property is equivalent to the 
capacity of $\partial\Omega$ relative to $\Omega$ being zero.
Alternatively these conditions are equivalent to the equality 
$S^D_t \varphi =S^N_t \varphi$ 
for all $\varphi \in L_2(\Omega)$ and $t > 0$,
where $S^N$ is the semigroup obtain by imposing Neumann boundary conditions 
on $\partial\Omega$.
This latter result is  related to the work of Arendt and Warma 
\cite{AW} \cite{AW2} on boundary conditions on the Laplacian on 
arbitrary domains and a number of our arguments are similar.

The analysis of \cite{RSi} was for a semigroup $S$ on $L_2(\Ri^d)$ 
generated by a second-order,
divergence-form, elliptic operator $H$ with $W^{1,\infty}$-coefficients and 
an open set $\Omega\subset \Ri^d$
with a Lipschitz boundary $\partial\Omega$.
Then it was established that $S_tL_2(\Omega)\subseteq L_2(\Omega)$ for all 
$t>0$ if and only if the capacity
$\capp(\partial\Omega)$ of $\partial\Omega$ measured with respect to the 
form $h$ associated with $H$ is zero.
It was also remarked that this equivalence fails if the coefficients of $H$ 
are not Lipschitz continuous.
The problem is that the degeneracy of the coefficients can differ depending 
whether one approaches
the boundary $\partial\Omega$ from $\Omega$ or from  $\overline\Omega^{\,\rm c}$.
The situation was clarified in \cite{ER29} by the demonstration that 
invariance could be completely
characterized by a condition on the  capacities relative to $\Omega$ and  
$\Omega^{\rm c}$ with no regularity 
required of  the coefficients
or  the boundary $\partial\Omega$.
In addition the set $\Omega$ is allowed to be measurable.
The results of \cite{ER29} were derived in the general framework of local 
Dirichlet forms and the 
current discussion will also be carried out in this framework.

We assume throughout  that $X$ is a locally compact $\sigma$-compact
metric space equipped with a positive Radon measure $\mu$ such that $\supp \mu = X$.
Let $\ce$ be a Dirichlet form on $X$.
The Dirichlet form is called {\bf regular} if $D(\ce) \cap C_c(X)$ is 
dense both in $D(\ce)$, with the graph norm, and in $C_0(X)$,
with the supremum norm.
Throughout this paper we  assume that 
$D(\ce) \cap C_c(X)$ is dense in $C_0(X)$.
Moreover, we also  require throughout that  $\ce$ is {\bf local} in the sense that 
$\ce(\psi,\varphi)=0$ for all 
$\varphi,\psi\in D(\ce)$ with $\varphi \, \psi=0$.
This notion appears slightly stronger  than locality as defined in \cite{FOT} 
but   if $\ce$ is regular then it is  equivalent  by a result of  Schmuland \cite{Schm}.
Let $\Omega$ be  an open subset of $X$.
We associate with the form $\ce$ a second form  $\ce_D $ which  
corresponds abstractly to $\ce$ with Dirichlet boundary conditions
imposed on $\partial\Omega$.
The latter form is defined by first setting
\[
D_\Omega = D_\Omega(\ce) 
= \{ \varphi\in D(\ce): \supp\varphi \mbox{ is a compact subset of } \Omega \}
\;\;\; .  \]
Since $D(\ce) \cap C_c(X)$ is dense in $C_0(X)$ it follows that $D_\Omega$ 
is dense in $L_2(\Omega)$.
Then we define $\overline{D_\Omega}$ as the closure of $D_\Omega$ 
with respect to the  graph norm on  $D(\ce)$.
Since $ \one_\Omega\, \varphi = \varphi$ for all $\varphi \in D_\Omega$ and the 
multiplication operator $\varphi \mapsto  \one_\Omega\, \varphi $ is continuous 
on $L_2(X)$,
it follows that $\overline{D_\Omega} \subseteq L_2(\Omega)$.
Here and in the sequel we identify $L_2(\Omega)$ in a natural way with the 
subspace $ \{  \one_\Omega\, \varphi  : \varphi \in L_2(X) \} $.
Finally $\ce_D (=\ce_{\Omega,  D})$ is defined as a form on $L_2(\Omega)$   
with domain $D(\ce_D) = \overline{D_\Omega}$
by $\ce_D = \ce|_{\overline{D_\Omega}}$.
Subsequently, in Section~\ref{Sord3}, we introduce a second form $\ce_N$ 
which corresponds to the introduction of 
Neumann boundary conditions on $\partial\Omega$.
But the definition of  $\ce_N$  is more complicated and its analysis 
requires stronger assumptions.
Therefore we first concentrate on the relatively simple form $\ce_D$.

It follows straightforwardly that  $\ce_D$ is a Dirichlet form on 
$L_2(\Omega)$ and $D(\ce_D) \cap C_c(\Omega)$ is dense in $C_0(\Omega)$.
Let $H_D (= H_{\Omega,D})$ and $S^D (=S^{\Omega,D})$ 
denote the operator and semigroup on 
$L_2(\Omega)$ associated with $\ce_D$. 
Since $S$ and $S^D$ are submarkovian semigroups they extend
to all the $L_p$-spaces including $L_\infty(X)$ and $L_\infty(\Omega)$.

Next we define the capacity and relative capacity of a set with respect
to the form $\ce$.
If $\Omega$ is a subset of $X$ and $A\subseteq \overline\Omega$ then 
the {\bf relative capacity} $\capp_\Omega(A) \in [0,\infty]$ is defined~by
\begin{eqnarray*}
\capp_\Omega(A)
= \capp_{\Omega,\ce}(A)
= \inf \{ \|\varphi\|_{D(\ce)}^2 & : & \varphi\in D(\ce)
    \mbox { and there exists an open } V \subset X  \nonumber\\
& & \mbox{such that }
    A \subset V \mbox{ and } \varphi\geq \one\mbox{ a.e.\ on } V \cap \Omega\}
\;\;\;.
\end{eqnarray*}
If $\Omega = X$ then $\capp(A) = \capp_\ce(A) = \capp_{X,\ce}(A)$ is the 
{\bf capacity} of the set $A$.
This version of relative capacity is the one used in \cite{ER29}, but it is 
probably different from the definition of relative capacity introduced earlier
by Arendt and Warma \cite{AW} \cite{AW2}.

If $\ce$ is regular and $\Omega$ is measurable 
then it follows from \cite{ER29}, Theorem~1.1, that $S$ leaves $L_2(\Omega)$ invariant
if and only if there exist $A_1,A_2\subseteq \partial\Omega$ such that 
$A_1\cup A_2=\partial\Omega$
and $\capp_\Omega(A_1)=0=\capp_{\Omega^{\rm c}}(A_2)$.
In particular, if  $\capp_\Omega(\partial\Omega)=0$
then $L_2(\Omega)$ is $S$-invariant.

Our main result gives a criterion for the validity of the converse 
of the latter statement.

\begin{thm}\label{tord340}
Adopt the foregoing definitions and assumptions.
Let $\Omega$ be an open subset of $X$.
Consider the following conditions.
\begin{tabel}
\item\label{tord340-1}
$S^{D}$ is conservative, i.e.\ $S^{D}_t\one_\Omega=\one_\Omega$ 
for all $t>0$.
\item\label{tord340-2}
$S_t\varphi =S^D_t\varphi$ for all  $\varphi\in L_2(\Omega)$ and $t>0$.
\item\label{tord340-3}
$\capp_\Omega(\partial\Omega)=0$.
\end{tabel}
Then {\rm \ref{tord340-1}$\Rightarrow$\ref{tord340-2}$\Rightarrow$\ref{tord340-3}}.
In particular Conditions~{\rm \ref{tord340-1}} and {\rm \ref{tord340-2}}
imply  that $L_2(\Omega)$ is $S$-invariant.

Moreover, if $S$ is conservative then {\rm \ref{tord340-2}$\Rightarrow$\ref{tord340-1}}.
Finally, if $\ce$ is regular, then {\rm \ref{tord340-3}$\Rightarrow$\ref{tord340-2}}.
\end{thm}

The theorem applies directly if $\ce$ is the form of a second-order, 
divergence-form, elliptic operator with real measurable
coefficients on $L_2(\Ri^d)$.
Then $\ce$  is regular, local and
the corresponding semigroup $S$ is conservative.
We will discuss this example more fully in Section~\ref{Sord4}.
The equivalence \ref{tord340-2}$\Leftrightarrow$\ref{tord340-3} 
generalizes a result of Arendt and Warma for the Laplacian 
(see \cite{AW2}, Proposition~2.5).

One can draw a stronger  conclusion if the capacity $\capp(\partial\Omega)=0$
and $\ce$ is regular, since this immediately implies that  
$\capp_\Omega(\partial\Omega)=0=\capp_{\overline\Omega^{\,\rm c}}(\partial\Omega)$.
There is, however, a converse to this statement if  $|\partial\Omega|=0$. 
Then the conditions 
$\capp_\Omega(\partial\Omega)=0=\capp_{\overline\Omega^{\,\rm c}}(\partial\Omega)$ 
imply that $\capp(\partial\Omega)=0$ by \cite{ER29}, Lemma~2.9.
(The condition $|\partial\Omega|=0$ is essential since 
$\capp(\partial\Omega)\geq |\partial\Omega|$.)
Therefore if $|\partial\Omega|=0$ then $\capp(\partial\Omega)=0$ is 
equivalent to both $S^D$ and $S^{\overline\Omega^{\,\rm c},D}$ being conservative
or to the conditions  $S_t\varphi =S^D_t\varphi$ and 
$S_t\psi= S^{\overline\Omega^{\,\rm c},D}_t\psi$ for all  
$\varphi\in L_2(\Omega)$, $\psi\in L_2(\overline\Omega^{\,\rm c})$ and $t>0$.

In Section~\ref{Sord3} we will give a further characterization of the condition $\capp_\Omega(\partial\Omega)=0$
in terms of Neumann boundary conditions.

\section{Dirichlet boundary conditions}\label{Sord2}

In this section we prove Theorem~\ref{tord340}.
The proof depends on  a couple of standard results which we use throughout this paper.

First, the $S$-invariance of $L_2(\Omega)$ is equivalent to the condition 
$\one_\Omega \, \varphi \in  D(\ce)$
for all $\varphi \in D(\ce)$ or for all $\varphi$ in a core of $\ce$.
These criteria are a corollary of a general result of Ouhabaz \cite{Ouh5}, Theorem~2.2, 
for local accretive forms 
(see also \cite{FOT}, Theorem~1.6.1, and \cite{ER29}, Proposition~2.1).
   
Secondly,  we need  an order relation between the semigroups $S$ and $S^D $.
Note that  each   bounded operator $A $ on $L_2(\Omega)$  can be extended to a 
bounded operator on $L_2(X)$,
still denoted by $A$, via 
$\varphi \mapsto A(\one_\Omega \, \varphi) \in L_2(\Omega) \subset L_2(X)$
for all $\varphi \in L_2(X)$.
In particular $S^D_t$ extends to a bounded operator on $L_2(X)$.
Note that  $\lim_{t\to0}S^D_t\varphi=\one_\Omega\varphi$ for all $\varphi \in L_2(X)$.

\begin{prop} \label{pord1}
If $\Omega$ is open
and $\varphi\in L_2(X)_+$
then $0 \leq S^D _t \varphi \leq S_t\varphi$
for all $t > 0$.
\end{prop}

The proposition follows from an adaptation of the reasoning of \cite{Are5}, Section~4.2.
Alternatively it can be deduced from \cite{Ouh5}, Theorem~2.24.
The proof relies on the following extension of Lemma~4.2.3 of \cite{Are5}.

\begin{lemma} \label{lord2}
Let $\varphi\in D(\ce_D )$ and $\psi\in D(\ce)_+$ satisfy
\begin{equation}
(\chi,\varphi) + \ce_D (\chi,\varphi)
\leq (\chi,\psi) +\ce(\chi,\psi)
\label{elord2;1}
\end{equation}
for all $\chi \in D(\ce_D )_+$.
Then $\varphi \leq \psi$.
\end{lemma}
\proof\
There exist $\varphi_1,\varphi_2,\ldots \in D_\Omega$ such that 
$\lim \|\varphi_n - \varphi\|_{D(\ce)} = 0$.
Then, however, $\supp (\varphi_n - \psi)_+ \subseteq \supp \varphi_n \subset \Omega$
since $\psi \geq 0$.
So $(\varphi_n - \psi)_+ \in D_\Omega$ for all $n \in \Ni$.
Moreover, $\lim (\varphi_n - \psi)_+ = (\varphi - \psi)_+$ in $D(\ce)$.
Hence $(\varphi - \psi)_+ \in D(\ce_D )$.

Secondly, set $\chi=(\varphi-\psi)_+$ in  (\ref{elord2;1}).
Then one deduces that 
\[
\|(\varphi-\psi)_+\|_2^2
=((\varphi-\psi)_+, \varphi-\psi)
\leq -\ce((\varphi-\psi)_+, \varphi-\psi)
=-\ce((\varphi-\psi)_+)
\leq 0
\;\;\;,
\]
where we used locality of $\ce$ in the last equality.
Hence $(\varphi-\psi)_+ = 0$ or, equivalently, $\varphi \leq \psi$.\hfill$\Box$

\ruimte

\noindent{\bf Proof of Proposition~\ref{pord1}\ }
Let $\tau \in L_2(\Omega)_+$.
Set $\varphi = (I+H_D )^{-1} \tau$ and $\psi = (I+H)^{-1} \tau$.
Then $\varphi\in D(H_D ) \subseteq D(\ce_D )$ and $\psi\in D(H)\subseteq D(\ce)$.
Moreover, $\psi\geq 0$ because $\tau\geq0$ and $S$ is submarkovian.
Now
\[
(\chi,\varphi) + \ce_D (\chi,\varphi)
= (\chi, (I + H_D ) \varphi)
= (\chi,\tau)
= (\chi,\psi) + \ce(\chi,\psi)
\]
for all $\chi\in D(\ce_D )$.
Therefore $(I+H_D )^{-1} \tau \leq (I+H)^{-1}\tau$
by Lemma~\ref{lord2}.
Similarly, $(I + \lambda H_D )^{-1} \tau \leq (I + \lambda H)^{-1} \tau$
for all $\lambda > 0$ and $\tau \in L_2(\Omega)_+$.
Then $S^D _t \tau \leq S_t \tau$ for all $t > 0$ since
$S_t \tau = \lim_{n \to \infty} (I + n^{-1} t H)^{-n} \tau$ with a similar expression for 
$S^D_t$.
Finally, since $S^D_t \tau = 0$ for all $\tau \in L_2(\Omega^{\rm c})$
the proposition follows.\hfill$\Box$

\begin{cor} \label{cord3}
If $\Omega_1 \subseteq \Omega_2$ are open
then $0 \leq S^{\Omega_1,D}_t \varphi \leq S^{\Omega_2,D}_t \varphi$
for all $\varphi\in L_2(X)_+$ and $t > 0$.
\end{cor}
\proof\
This follows from Proposition~\ref{pord1} with  $X$ replaced by $\Omega_2$,
$\ce$ replaced  by $\ce_{\Omega_2}$ and $S$ by $ S^{\Omega_2,D}$.\hfill$\Box$

\ruimte

Now we turn to the proof of Theorem~\ref{tord340}.

\ruimte

\noindent{\bf Proof of Theorem~\ref{tord340}\ }
``\ref{tord340-1}$\Rightarrow$\ref{tord340-2}''.
Let $\varphi \in L_1(\Omega) \cap L_2(\Omega)_+$ and $t > 0$.
Then $S^D _t \varphi \leq S_t \varphi$ by Proposition~\ref{pord1}.
Therefore using Condition~\ref{tord340-1} and the positivity and contractivity of $S$ one has
\begin{eqnarray}
\|\varphi\|_1
= (\one_\Omega, \varphi)
= (S^D _t \one_\Omega, \varphi)
= (\one_\Omega, S^D _t \varphi)
 \leq  (\one_\Omega, S_t \varphi)
\leq  (\one, S_t \varphi)  
 =  \|S_t \varphi\|_1
\leq \|\varphi\|_1
\; . \;\;\;\label{ecord3.1;2}\end{eqnarray}
Hence all three inequalities are in fact equalities.
Since the second inequality in (\ref{ecord3.1;2}) is an equality it follows that 
$(\one_{\Omega^{\rm c}}, S_t \varphi) = 0$.
Therefore $\one_{\Omega^{\rm c}} S_t \varphi = 0$ and 
$S_t \varphi \in L_2(\Omega)$.
Since the first inequality in (\ref{ecord3.1;2}) is an equality one deduces 
from the order relation  $S^D _t \varphi \leq S_t \varphi$ of Proposition~\ref{pord1}
that $S^D _t \varphi = S_t \varphi$.
But this immediately implies that $S^D _t \psi = S_t \psi$ for all $t > 0$ and 
$\psi \in L_2(\Omega)$.
Thus Condition~\ref{tord340-2} is established.

``\ref{tord340-2}$\Rightarrow$\ref{tord340-3}''.
If $\varphi \in L_2(\Omega)$ then $S_t \varphi = S^D _t \varphi \in L_2(\Omega)$.
So $L_2(\Omega)$ is $S$-invariant.
Next, let $K \subset X$ compact.
Since $X$ is locally compact and $D(\ce) \cap C_c(X)$ is dense in $C_0(X)$ there 
exist an open set $V$ and a $\varphi \in D(\ce) \cap C_c(X)$
such that $\varphi \geq \one_V \geq \one_K$ pointwise.
Then $ \one_\Omega\,\varphi \in D(\ce) \cap L_2(\Omega)$, by $S$-invariance of $L_2(\Omega)$,  and 
\[
\lim_{t \downarrow 0} t^{-1} (\one_\Omega\,\varphi , (I - S^D _t) \one_\Omega\,\varphi )
= \lim_{t \downarrow 0} t^{-1} (\one_\Omega\,\varphi , (I - S_t) \one_\Omega\,\varphi )
= \ce(\one_\Omega\,\varphi )
\]
exists.
Therefore $\one_\Omega\,\varphi \in D(\ce_D )$.
By definition of $\ce_D $ there exist $\psi_1,\psi_2,\ldots \in D_\Omega$ such that 
$\lim_{n \to \infty} \psi_n =\one_\Omega\,\varphi $ in $D(\ce)$.
Then $\one_\Omega\,\varphi  - \psi_n \in D(\ce)$, 
$K \cap \partial \Omega \subset V \setminus \supp \psi_n$,
the set $V \setminus \supp \psi_n$ is open and 
$\one_\Omega\,\varphi  - \psi_n \geq \one$ a.e.\ on $(V \setminus \supp \psi_n) \cap \Omega$
for all $n \in \Ni$.
Therefore $\capp_\Omega(K \cap \partial\Omega) \leq \|\one_\Omega\,\varphi  - \psi_n\|_{D(\ce)}^2$
for all $n \in \Ni$ and $\capp_\Omega(K \cap \partial\Omega) = 0$.
Since $X$ is $\sigma$-compact one deduces that $\capp_\Omega(\partial\Omega) = 0$.

``\ref{tord340-2}$\Rightarrow$\ref{tord340-1}''.
Suppose that $S$ is conservative.
If $\varphi \in L_1(\Omega) \cap L_2(\Omega)$ then
\[
(\varphi, S^D _t \one_\Omega)
= (S^D _t \varphi, \one_\Omega)
= (S_t \varphi, \one_\Omega)
= (\varphi, S_t \one_\Omega)
= (\varphi, \one_\Omega S_t \one)
= (\varphi, \one_\Omega)
\]
for all $t > 0$.
Therefore $S^D _t \one_\Omega = \one_\Omega$ for all $t > 0$.

``\ref{tord340-3}$\Rightarrow$\ref{tord340-2}''.
Finally, suppose that $\ce$ is regular.
We shall prove that if $\varphi\in D(\ce)$ then  $\one_\Omega\,\varphi  \in D(\ce_D )$.
We argue as in the proof of \cite{ER29}, Theorem~2.4.

Since $\capp_\Omega(\partial\Omega)=0$ 
for all $n\in \Ni$ there exist $\psi_n\in D(\ce)$
and an open $V_n\subset X$ such that $\partial\Omega\subset V_n$, $\psi_n\geq \one$ almost everywhere
on $V_n\cap \Omega$ and $\|\psi_n\|_{D(\ce)}\leq 1/n$.
Without loss of generality we may assume that $0 \leq \psi_n \leq \one$.
Let $\varphi \in D(\ce) \cap C_c(X)$.
Let $n \in \Ni$.
Define $\varphi_n = (\varphi - \varphi \, \psi_n) \one_\Omega \in L_2(\Omega)$.
Then $\supp \varphi_n$ is compact and 
\[
\supp \varphi_n
\subset \overline{\Omega \cap V_n^{\rm c}} 
\subset \overline \Omega \cap V_n^{\rm c}
\subset \Omega
\;\;\; .  \]
Hence there exists a $\chi \in D(\ce) \cap C_c(\Omega)$ such that 
$\chi|_{\supp \varphi_n} = \one$.
Then $\varphi_n = (\varphi - \varphi \, \psi_n) \chi \in D(\ce)$.
So $\varphi_n\in D_\Omega\subseteq D(\ce_D )$.
It follows from locality that 
\begin{eqnarray*}
\ce(\varphi_n)
& \leq & \ce(\varphi_n) + \ce((\varphi - \varphi \, \psi_n) \one_{\Omega^{\rm c}})  \\
& = & \ce(\varphi - \varphi \, \psi_n)
\leq 2\, \ce(\varphi) + 4\, \ce(\varphi) \, \|\psi_n\|_\infty^2 
    + 4 \,\ce(\psi_n) \, \|\varphi\|_\infty^2
\leq 6\, \ce(\varphi) + 4\, \|\varphi\|_\infty^2
\end{eqnarray*}
for all $n\in\Ni$.
So the sequence $\varphi_1,\varphi_2,\ldots$ has a weakly convergent subsequence 
$\varphi_{n_1},\varphi_{n_2},\ldots$ in the Hilbert space $D(\ce_D)$.
Clearly $\lim_{n \to \infty} \varphi_n = \one_\Omega \, \varphi$ in $L_2(\Omega)$.
So $\one_\Omega \, \varphi \in D(\ce_D)$ for all $\varphi \in D(\ce) \cap C_c(X)$.

Since $D(\ce) \cap C_c(X)$ is dense in $D(\ce)$ by regularity it follows that 
$L_2(\Omega)$ is $S$-invariant by \cite{ER29}, Proposition~2.1.
Moreover, by density,
$D(\ce) \cap L_2(\Omega) \subseteq D(\ce_D)$.
Since the converse inclusion is obvious it follows that  $D(\ce)\cap L_2(\Omega)= D(\ce_D)$.
Hence  $S_t\varphi=S^D_t\varphi$ for all $\varphi\in L_2(\Omega)$ and $t>0$.
This completes the proof of Theorem~\ref{tord340}.\hfill$\Box$

\section{Neumann boundary conditions}\label{Sord3}

The  form corresponding to $\ce$ with  Neumann boundary conditions
on $\partial\Omega$ is defined in terms  of the truncations of $\ce$.
If  $\chi \in D(\ce) \cap L_\infty(X)_+$ 
then the {\bf truncated form} $\ce_\chi$ is given by
$D(\ce_\chi)= D(\ce) \cap L_\infty(X )$ and 
\[
\ce_\chi(\varphi) = \ce(\chi\,\varphi, \varphi) - 2^{-1} \ce(\chi,\varphi^2)
\] 
for all $\varphi\in D(\ce) \cap L_\infty(X)$.
It has three basic properties:
\begin{equation}
0\leq \ce_\chi(\varphi) \leq \|\chi\|_\infty \, \ce(\varphi)
\;\;\;,
\label{eord3.2}
\end{equation}
\begin{equation}
\ce_\chi(0\vee\varphi\wedge \one) \leq \ce_\chi(\varphi)
\;\;\;,
\label{eord3.3}
\end{equation}
and
\begin{equation}
\mbox{ if  }  0\leq \chi_1\leq \chi_2 \mbox{ then }0\leq \ce_{\chi_1}(\varphi)\leq \ce_{\chi_2}(\varphi)
\label{eord3.4}
 \end{equation}
 where all three properties are valid for all $\varphi\in D(\ce) \cap L_\infty(X)$.
These properties are established in  \cite{BH}, Proposition~I.4.1.1.

It follows from (\ref{eord3.2}) that $\ce_\chi$ can be extended to $D(\ce)$ by continuity.
The extension, which we continue to denote by $\ce_\chi$, still satisfies the Markovian property (\ref{eord3.3})
and the monotonicity property (\ref{eord3.4}).

Next  for each open subset  $\Omega$   of $X$  define the convex subset $\cc_\Omega$ of $D(\ce)$ by
\[
\cc_\Omega
=\{\chi\in D(\ce) \cap L_\infty(X )\,,\,0\leq\chi\leq \one_\Omega\}
\;\;\;.
\]
It follows that $\cc_\Omega$ is a directed set with respect to the natural order. 
In particular if $\chi_1, \chi_2\in \cc_\Omega$ then $\chi_{12}=\chi_1+\chi_2-\chi_1\chi_2\in \cc_\Omega$.
Moreover, $\chi_{12}-\chi_1=\chi_2(\one_\Omega-\chi_1)\geq0$ and $\chi_{12}-\chi_2=\chi_1(\one_\Omega-\chi_2)\geq0$.
Therefore it follows from (\ref{eord3.4}) that 
$\chi \mapsto \ce_\chi$ is a monotonically increasing net of quadratic forms 
with the common domain $D(\ce)$.
Then one can define a  form $\ce_N\,(=\ce_{\Omega, N})$ by $D(\ce_N)= D(\ce) $ and 
\[
\ce_N(\varphi)
=\lim_{\chi\in\cc_\Omega}\ce_{ \chi}(\varphi)
=\sup \{ \ce_{ \chi}(\varphi) : \chi \in \cc_\Omega \}
\;\;\;.
\]
Since $\ce_N$ is defined as a limit of quadratic forms it is automatically a quadratic form on $L_2(X)$ and 
it follows from (\ref{eord3.2}) and (\ref{eord3.3}) 
that $\ce_N$ satisfies the continuity property
\begin{equation}
0\leq \ce_N(\varphi)\leq \ce(\varphi)
\label{eord3.5}
\end{equation}
and the Markovian property
\begin{equation}
\ce_N(0\vee\varphi\wedge \one)\leq \ce_N(\varphi)
\label{eord3.6}
\end{equation}
for all $\varphi\in D(\ce)$.
We emphasize that $\ce_N$ is a form on $L_2(X)$.

The definition of $\ce_N$ is motivated by the theory of second-order elliptic operators.
Let $X=\Ri$ and define $\ce$ by  $D(\ce)=W^{1,2}(\Ri)$ and $\ce(\varphi)=\int_\Ri|\varphi'|^2$.
Then $\ce_\chi(\varphi)=\int_\Ri\chi|\varphi'|^2$ and $\ce_N(\varphi)=\int_\Omega|\varphi'|^2$.

Our aim is to compare the forms $\ce_D$ and $\ce_N$ on $L_2(\Omega)$ but in general
$\ce_N$  is  not closed nor even closable.
In fact it is closed under quite general assumptions (see Proposition~\ref{pord3.2} below) but 
in any case one can introduce the   relaxation $\widehat{\ce_N}$ of $\ce_N$.

The {\bf relaxation}
$\hat t$ of a quadratic form $t$ is variously called the  lower semi-continuous 
regularization (see \cite{ET}, page~10) or the  relaxed form 
(see \cite{DalM}, page~28).
It is the closure of the largest closable form which is less than or equal to $t$
(see  \cite{bSim5} Theorem~2.2).
In particular, if $t$ is closable then $\hat t$ is the closure.

The relaxation  $\widehat{\ce_N}$ of $\ce_N$ is automatically a Dirichlet form;
it is positive, closed and   satisfies (\ref{eord3.6}).
Moreover, it  satisfies $\widehat{\ce_N}(\varphi)\leq \ce(\varphi)$ for all 
$\varphi\in D(\ce)$ by (\ref{eord3.5}).
Let $H_N\,(=H_{\Omega,N})$ and $S^N\,(=S^{\Omega,N})$ denote the operator and 
submarkovian semigroup on $L_2(X)$
associated with~$\widehat{\ce_N}$.

\begin{remarkn} \label{rord355}
If $\varphi \in D(\ce) \cap C_c(\overline \Omega^{\,\rm c})$ then 
$\ce_\chi(\varphi) = 0$ for all $\chi \in \cc_\Omega$ by locality.
Therefore $\ce_N(\varphi) = 0$ and $\widehat{\ce_N}(\varphi) = 0$.
Since $\varphi \in D(\ce) \cap C_c(\overline \Omega^{\,\rm c})$ is dense in 
$C_c(\overline \Omega^{\,\rm c})$ and $C_c(\overline \Omega^{\,\rm c})$ is dense in 
$L_2(\overline \Omega^{\,\rm c})$ one deduces that $\widehat{\ce_N}(\varphi) = 0$
for all $\varphi \in L_2(\overline \Omega^{\,\rm c})$. 
Hence if $\varphi \in D(\widehat{\ce_N})$ then $\one_\Omega \varphi \in D(\widehat{\ce_N})$
and $\widehat{\ce_N}(\varphi) = \widehat{\ce_N}(\one_\Omega \varphi)$.
In particular, the space $L_2(\overline \Omega)$ is invariant under~$S^N$.
\end{remarkn}

\begin{prop}\label{pord3.10}
Let $\Omega \subset \Ri^d$ be open.
If  $S^{D}_t\varphi=S^N_t\varphi$ for all $\varphi\in L_2(\Omega)$ and $t>0$
then $L_2(\Omega)$ is $S^N$-invariant and $\capp_\Omega(\partial\Omega)=0$.
\end{prop}
\proof\
Since $S^D$ leaves $L_2(\Omega)$ invariant the $S^N$-invariance follows immediately.
But the latter property implies that if $\varphi\in D(\widehat{\ce_N})$ then 
$\one_\Omega\varphi\in D(\widehat{\ce_N})$.
Next let $\varphi\in D(\ce)$.
Then $\varphi \in D(\widehat{\ce_N})$ and 
\[
\lim_{t \downarrow 0} t^{-1} (\one_\Omega\,\varphi , (I - S^D _t) \one_\Omega\,\varphi )
= \lim_{t \downarrow 0} t^{-1} (\one_\Omega\,\varphi , (I - S^N_t) \one_\Omega\,\varphi )
= \widehat{\ce_N}(\one_\Omega\,\varphi )<\infty
\;\;\; .  \]
So $\one_\Omega\,\varphi \in D(\ce_D )$.
The rest of the proof is then a repetition of the argument that 
\ref{tord340-2}$\Rightarrow$\ref{tord340-3} 
in Theorem~\ref{tord340}.\hfill$\Box$

\ruimte

Under more stringent assumptions (see Theorem~\ref{tord341}) we will prove that  
Proposition~\ref{pord3.10} has a converse.
One key condition is strong locality.

We define   $\ce$ to be {\bf strongly local} if $\ce(\varphi,\psi)=0$ for all 
$\varphi,\psi\in D(\ce)$
and $a\in\Ri$ such that $(\varphi+a\one)\psi=0$.
This condition corresponds to locality  in the sense of \cite{BH}.

Strong locality gives a couple of useful implications.

\begin{lemma} \label{lord350}
Suppose $\ce$ is strongly local and regular.
Then 
\[
\ce_N(\varphi)
= \sup \{ \ce_\chi(\varphi) : \chi \in D(\ce) \cap C_c(\Omega), \; 0 \leq \chi \leq \one \}
\]
for all $\varphi \in D(\ce)$.
\end{lemma}
\proof\
First notice that there are 
$\chi_1,\chi_2,\ldots \in D(\ce) \cap C_c(\Omega)_+$ such that 
$\chi_n \uparrow \one_\Omega$.
Then $\ce_N(\varphi) = \lim_{n \to \infty} \ce_{\chi_n}(\varphi)$
for all $\varphi \in D(\ce)$
(see the discussion on page~82 in \cite{ERS4}, which requires $\ce$
to be regular and strongly local).\hfill$\Box$

\ruimte

Thus if $\ce$ is regular and strongly local then 
one can replace the set $\cc_\Omega$ by the 
set $ \{ \chi \in D(\ce) \cap C_c(\Omega) : 0 \leq \chi \leq \one \} $
in the definition of $\ce_N$.

Next we establish that if   $\ce$ is strongly local then there is an order relation 
between $S^D$ and $S^N$.

\begin{prop} \label{pord340}
Let $\Omega \subset \Ri^d$ be open.
If $\ce$ is strongly local then $\ce_D \subseteq \ce_N$.
Moreover, 
$0 \leq S^D_t \varphi \leq S^N_t \varphi$ for all $\varphi \in L_2(X)_+$ 
and $t > 0$.
\end{prop}
\proof\
Clearly $\ce_N(\varphi) = \ce(\varphi)$ for all $\varphi \in D_\Omega$.
But $D_\Omega$ is dense in $D(\ce_D)$.
Hence it follows from (\ref{eord3.5}) that 
$\ce_N(\varphi) = \ce(\varphi)=\ce_D(\varphi)$ for all $\varphi \in D(\ce_D)$.
Therefore  $\ce_D \subseteq \ce_N$.

Since  $L_2(\Omega)$ is $S^D$-invariant and $S^D_t \varphi = 0$ for all
$\varphi \in L_2(\Omega^{\rm c})$, by definition, 
it suffices to prove the order property of the semigroups
for all $\varphi\in L_2(\Omega)_+$.

Let $\varepsilon > 0$ and define the form $\ce_{N \varepsilon}$ by 
$\ce_{N \varepsilon} = \ce_N + \varepsilon \, \ce$.
Then $\ce_{N \varepsilon}$ is a Dirichlet form and $D(\ce_{N \varepsilon}) \cap C_c(X)$
is dense in $C_0(X)$.
Moreover, $D_\Omega(\ce) = D_\Omega(\ce_{N \varepsilon})$
and $(\ce_{N \varepsilon})_D = (1 + \varepsilon) \ce_D$.
Therefore we can apply Proposition~\ref{pord1} to deduce that 
$0 \leq S^D_{(1+\varepsilon)t} \varphi \leq S^{N \varepsilon}_t \varphi$
for all $t > 0$ and $\varphi \in L_2(\Omega)_+$, where 
$S^{N \varepsilon}$ is the semigroup associated with the Dirichlet form
$\ce_{N \varepsilon}$.
Since $\lim_{\varepsilon \downarrow 0} S^{N \varepsilon}_t = S^N_t$ strongly 
for all $t > 0$ by \cite{Kat1}, Theorem~VIII.3.11, the proposition is established.\hfill$\Box$

\begin{remarkn} \label{rord3}
The semigroup domination property of Proposition~\ref{pord340} can be characterized in 
terms of the forms
$\ce_D$ and $\widehat{\ce_N}$ by a general result of Ouhabaz (see \cite{Ouh5}, Theorem~2.24).
In particular it follows that $D(\ce_D)$ is an ideal of $D(\widehat{\ce_N})$, i.e.\ if 
$0\leq \varphi\leq \psi$
with $\varphi\in D(\widehat{\ce_N})$ and $\psi\in D(\ce_D)$ then $\varphi\in D(\ce_D)$.
\end{remarkn}

Under the additional assumption that $\ce$ is regular one can deduce that $S$-invariance of 
$L_2(\Omega)$ suffices for  equality of $S$ and $S^N$ in restriction to  $L_2(\Omega)$.

\begin{prop}\label{pord3.2} 
Let $\Omega \subset \Ri^d$ be open.
Assume $\ce$ is regular and strongly local and that $L_2(\Omega)$ is $S$-invariant.
Then
\[
\ce_N(\varphi)=\ce(\one_\Omega\varphi)
\]
for all $\varphi\in D(\ce)$.
Therefore $ \ce_N$ is closed.
Moreover, $S^N$ leaves $L_2(\Omega)$ invariant and 
$S^N_t\varphi=S_t\varphi$ for all $\varphi\in L_2(\Omega)$ and $t>0$.
\end{prop}
\proof\
Fix $\varphi\in D(\ce)\cap C_c(X)$.
Then $\one_\Omega\varphi\in D(\ce)$ since $L_2(\Omega)$ is $S$-invariant.
Moreover, $\chi=\one_\Omega\chi$ for all $\chi\in \cc_\Omega$.
Therefore
\begin{eqnarray*}
\ce_N(\varphi)&=&\lim_{\chi\in \cc_\Omega}\Big(\ce(\varphi,\chi\varphi)-2^{-1}\ce(\chi,\varphi^2)\Big)\\[5pt]
&=&\lim_{\chi\in \cc_\Omega}\Big(\ce(\varphi,\chi\one_\Omega\varphi)-2^{-1}\ce(\one_\Omega\chi,\varphi^2)\Big)\\[5pt]
&=&\lim_{\chi\in \cc_\Omega}\Big(\ce(\one_\Omega\varphi,\chi\one_\Omega\varphi)-2^{-1}\ce(\chi,(\one_\Omega\varphi)^2)\Big)
=\ce_N(\one_\Omega\varphi)
\end{eqnarray*}
where we have used  locality.
Thus $\ce_N(\varphi)=\ce_N(\one_\Omega\varphi)$.

Next choose $\psi\in D(\ce)\cap C_c(X)$ with $\psi\geq \one_K$ where  $K=\supp\varphi$.
Then, replacing $\psi$ by $0\vee \psi\wedge \one$ if necessary, one can assume 
$0\leq \psi\leq \one$ and $\psi=\one$ on $K$.
Set $\chi=\one_\Omega\psi$.
Then $\chi\in\cc_\Omega$ and $\chi=\one$ on $K\cap\Omega$.
Therefore  $\chi\one_\Omega\varphi=\one_\Omega\varphi$ and 
\[
\ce_\chi(\one_\Omega\varphi)=\ce(\one_\Omega\varphi,\chi\one_\Omega\varphi)-2^{-1}\ce(\chi,(\one_\Omega\varphi)^2)
=\ce(\one_\Omega\varphi)
\]
 by strong locality.
Hence $\ce_N(\one_\Omega\varphi)=\ce(\one_\Omega\varphi)$.
But in the previous paragraph we established that 
$\ce_N(\varphi)=\ce_N(\one_\Omega\varphi)$.
Therefore  $\ce_N(\varphi)=\ce(\one_\Omega\varphi)$ for all $\varphi\in D(\ce)\cap C_c(X)$.
This equality then extends to all $\varphi\in D(\ce)$ by regularity of $\ce$.
The remaining statements of the Proposition~\ref{pord3.2} are straightforward.\hfill$\Box$

\ruimte

We now prove a kind of converse of Proposition~\ref{pord3.10}.

\begin{thm}\label{tord341} 
Assume $\ce$ is regular and   strongly local.
The following conditions are equivalent.
\begin{tabel}
\item\label{tord341-2}
$\capp_\Omega(\partial\Omega)=0$.
\item\label{tord341-1}
$S^{D}_t\varphi=S^N_t\varphi$ for all $\varphi\in L_2(\Omega)$ and all  $t>0$.
\end{tabel}
\end{thm}
\proof\
The implication \ref{tord341-1}$\Rightarrow$\ref{tord341-2} is established 
by Proposition~\ref{pord3.10} without the regularity and strong locality.

``\ref{tord341-2}$\Rightarrow$\ref{tord341-1}''.
Suppose $\capp_\Omega(\partial\Omega)=0$.
Then the implication
{\rm \ref{tord340-3}$\Rightarrow$\ref{tord340-2}} of Theorem~\ref{tord340} gives
$S^{D}_t\varphi=S_t\varphi $ for all  $\varphi\in L_2(\Omega)$ and $t>0$.
But $\capp_\Omega(\partial\Omega)=0$ also implies that $L_2(\Omega)$ is $S$-invariant.
Therefore  Proposition~\ref{pord3.2} 
gives $S_t\varphi =S^N_t\varphi$  for all  $\varphi\in L_2(\Omega)$ and $t>0$.
Hence by combination of these conclusions one obtains Statement~\ref{tord341-1} of the 
theorem.\hfill$\Box$

\ruimte

It follows from Theorems~\ref{tord340} and \ref{tord341} that the 
relative capacity condition $\capp_\Omega(\partial\Omega)=0$
is equivalent to $S^{D}_t\varphi=S_t\varphi$ for all 
$\varphi\in L_2(\Omega)$ and $t>0$ or to
$S^{D}_t\varphi=S^N_t\varphi$ for all $\varphi\in L_2(\Omega)$ and $t>0$.
It is not equivalent, however, to  $S^N_t\varphi=S_t\varphi$ 
for all $\varphi\in L_2(\Omega)$ and $t>0$.
A counterexample can be given as follows.
Define the form  $h$ on $L_2(\Ri)$ by 
$D(h)=L_2(-\infty,0)\oplus W^{1,2}(0,\infty)$ and $h(\varphi)=\int^\infty_0|\varphi'|^2$.
Set $\Omega=\langle0,\infty\rangle$.
Then $h=h_N$ and $S=S^N$.
But $h_D$ is the restriction of $h$ to $W^{1,2}_0(0,\infty)$.
Therefore $S^D$ is not the restriction of $S$ to $L_2(\Omega)$ and $\capp_\Omega(\{0\})\neq0$.
In fact $\capp_\Omega(\{0\})\geq (4\pi)^{-1}$.

\section{Degenerate elliptic  operators} \label{Sord4}

The foregoing results can be applied to degenerate elliptic operators
on $\Ri^d$.

Let $(c_{kl})$ be a symmetric $d\times d$-matrix with coefficients 
$c_{kl} \in L_\infty(\Ri^d)$ 
such that  $C(x) = (c_{kl}(x))$ is positive-definite for almost all $x \in \Ri^d$.
Define   the  positive quadratic form $h$ by $D(h) = W^{1,2}(\Ri^d)$ and
\[
h(\varphi)
= \sum_{k,l=1}^d (\partial_k \varphi, c_{kl} \, \partial_l \varphi)
\;\;\; .  \]
We call $h$ the degenerate elliptic form with coefficients $(c_{kl})$.
Further let $\hat h$ denote the relaxation of $h$.
It is established in  \cite{ERSZ2}, Theorem~1.1, that 
$\hat h$ is a regular, strongly local, Dirichlet form.
(The relaxation is referred to as the viscosity form in \cite{ERSZ2} and the definition of
locality used in  this reference corresponds to strong locality as defined in Section~\ref{Sord3}.)
Moreover, the submarkovian semigroup $S$ associated with $\hat  h$ is conservative 
by Theorem~3.7 of \cite{ERSZ1}.
Therefore all the statements of Theorems~\ref{tord340} and \ref{tord341} are equivalent 
for  $\hat h$ and  the corresponding elliptic operator $H$
and submarkovian semigroup~$S$.

The form $(\hat h)_D$ corresponding to $\hat h$ and  the open subset
$\Omega\subseteq \Ri^d$  is the $D(\hat h)$-closure  of the restriction of $\hat h$ 
to $C_c^\infty(\Omega)$.
Therefore, if $c_{kl} = \delta_{kl}$, i.e.\ if $h\,(=\hat h)$ is the form of the Laplacian,
then $h_D$ corresponds to the usual Laplacian with Dirichlet boundary conditions.
It is not the case, however, that $\widehat{h_N}$  always corresponds to the Laplacian  
with Neumann boundary conditions.
The form $l_N$ of the Laplacian with Neumann boundary conditions
is usually defined with the domain $D(l_N)=W^{1,2}(\Omega)$.
Note that $D(l_N) \subset L_2(\Omega)$.
But, by Remark~\ref{rord355},  one can always  write $D(\widehat{h_N}) = D_{\Omega,N} \oplus L_2(\overline \Omega^{\,\rm c})$
with $D_{\Omega,N}$ a subspace of $L_2(\overline \Omega)$.
If $|\partial \Omega| > 0$ then clearly $D_{\Omega,N} \neq W^{1,2}(\Omega) = D(l_N)$
and $D_{\Omega,N}$ contains elements which are not in $D(l_N)$.
But $D_{\Omega,N}$ can be a 
strict subset of $W^{1,2}(\Omega)$ even if $|\partial \Omega| = 0$.
If, for example, $d=1$ and $\Omega=\langle-1,0\rangle\cup \langle 0,1\rangle$ then 
$D_{\Omega,N} = W^{1,2}(-1,1) \subsetneqq W^{1,2}(-1,0) \oplus W^{1,2}(0,1) 
= W^{1,2}(\Omega)$.

Theorem~\ref{tord340} gives, in principle, a practical way of concluding that the 
semigroup $S^D$ corresponding to Dirichlet boundary conditions is conservative 
on $L_\infty(\Omega)$.
It suffices to verify that $\capp_\Omega(\partial\Omega)=0$.
But calculating the relative capacity is not straightforward.
The next proposition gives sufficient and practical 
conditions to make the verification.

\begin{prop} \label{pord421}
Let $\Omega \subset \Ri^d$ be open.
Let $h_1$ and $h_2$ be degenerate elliptic forms with coefficients 
$(c^{(1)}_{kl})$ and $(c^{(2)}_{kl})$.
Suppose there exists an $a \in \Ri$ such that 
$C^{(1)}(x) \leq a \, C^{(2)}(x)$ for almost every $x \in \Omega$.
Moreover, suppose that $\capp_{\Omega,\widehat{h_2}}(\partial \Omega) = 0$.
Then the semigroup $S^{(1)D}$ associated with $(\widehat{h_1})_D$ is conservative.
\end{prop}

The proof relies on the fact that the relaxation depends locally on the 
coefficients of the form.

\begin{lemma} \label{lord420}
Let $h_1$ and $h_2$ be degenerate elliptic forms with coefficients 
$(c^{(1)}_{kl})$ and $(c^{(2)}_{kl})$.
Let $U \subset \Ri^d$ be an open set and suppose that 
$c^{(1)}_{kl}|_U = c^{(2)}_{kl}|_U$ for all $k,l \in \{ 1,\ldots,d \} $.
Let $\varphi \in L_2(\Ri^d)$ and suppose that $\supp \varphi \subset U$.
Then $\varphi \in D(\widehat{h_1})$ if and only if $\varphi \in D(\widehat{h_2})$
and in this case $\widehat{h_1}(\varphi) = \widehat{h_2}(\varphi)$.
\end{lemma}
\proof\
Without loss of generality we may assume that $C^{(2)}(x) = 0$
for all $x \in U^{\rm c}$.
Vogt \cite{Vog1} proved that there exists a measurable function
$p \colon \Ri^d \to \Ri^{d \times d}$, with values in the 
orthogonal projections, such that the degenerate elliptic form $k_1$
with coefficients
$x \mapsto (C^{(1)}(x))^{1/2} \, p(x) \, (C^{(1)}(x))^{1/2}$ is closable
and $\widehat{h_1} = \overline{k_1}$.
Following the constructive proof  in \cite{Vog1} it follows that 
the degenerate elliptic form $k_2$
with coefficients
$x \mapsto (C^{(1)}(x))^{1/2} \, \one_U(x) \, p(x) \, (C^{(1)}(x))^{1/2}$ is closable
and $\widehat{h_2} = \overline{k_2}$.
Then the rest of the proof of the lemma is clear.\hfill$\Box$

\ruimte

\noindent
{\bf Proof of Proposition~\ref{pord421}}\hspace{5pt}\
For all $k,l \in \{ 1,\ldots,d \} $ define $c^{(3)}_{kl} \colon \Ri^d \to \Ri$
by $c^{(3)}_{kl} = \one_\Omega \, c^{(1)}_{kl}$.
Then $C^{(3)}(x) \leq a \, C^{(2)}(x)$ for almost every $x \in \Ri^d$.
Hence $\capp_{\Omega, \widehat{h_3}}(\partial \Omega) = 0$, where
$h_3$ is the degenerate elliptic forms with coefficients 
$(c^{(3)}_{kl})$.
Therefore the semigroup $S^{(3)D}$ associated with $(\widehat{h_3})_D$ is conservative.
But $(\widehat{h_3})_D = (\widehat{h_1})_D$ by Lemma~\ref{lord420}.
Hence the semigroup $S^{(1)D}$ associated with $(\widehat{h_1})_D$ is 
conservative.\hfill$\Box$

\ruimte

The assumptions of Proposition~\ref{pord421} are satisfied in many 
cases, see \cite{ER29} Section~3, or under the more stringent condition
$\capp_{\widehat{h_2}}(\partial \Omega) = 0$ see \cite{RSi}.

\subsection*{Acknowledgement}
The authors are  indebted to Adam Sikora for explaining  that
 Dirichlet boundary conditions could be used to characterize  subspaces 
invariant under the semigroup generated
by a degenerate elliptic operator and for sketching a proof  based on
finite propagation speed methods.

Part of the work was carried out whilst the first author was visiting 
the Australian National University
with partial support from the Centre for Mathematics and it 
Applications and part of the work was
carried out whilst the second author was visiting the University of Auckland 
with financial support
from the Faculty of Science.

\end{document}